\documentclass[reqno]{amsart}

\usepackage{enumerate}
\usepackage{ifpdf}
\usepackage{amsmath}
\usepackage{amsfonts}
\usepackage{amssymb}
\usepackage{amsthm}
\usepackage[hyperfootnotes=false]{hyperref}
\usepackage{setspace}
\usepackage{amsrefs}

\newcommand{\ignore}[1]{}

\newcommand{\C}{{\mathbb{C}}}
\newcommand{\R}{{\mathbb{R}}}

\newcommand{\sF}{{\mathcal{F}}}

\newcommand{\sM}{{\mathcal{M}}}
\newcommand{\sO}{{\mathcal{O}}}

\newcommand{\rank}{\operatorname{rank}}

\newtheorem{thm}{Theorem}[section]
\newtheorem*{thmnonum}{Theorem}

\newtheorem{lemma}[thm]{Lemma}

\newtheorem{question}[thm]{Question}

\theoremstyle{definition}

\newtheorem{example}[thm]{Example}

\theoremstyle{remark}

\author{Ji\v{r}\'i Lebl}
\address{Department of Mathematics, University of Illinois
at Urbana-Champaign, 
Urbana, IL 61801, USA}
\email{jlebl@math.uiuc.edu}
\date{\today}

\ifpdf
\fi

\title
{Pullback of varieties by finite maps}

\begin{document}


\begin{abstract}
We study the local geometry of the pullback of a variety via a
finite holomorphic map.  In particular, we are looking for properties of
$V = F^{-1}(W)$ such that if $V$ has the property $A$, then
$W$ must have the property $A$.  We show that $A$ can be the property of
normality or prefactoriality.  We also show that $A$ can be the property
of smoothness, under extra assumptions.
\end{abstract}

\maketitle



\section{Introduction} \label{section:intro}

When a complex analytic variety
is pulled back by an arbitrary holomorphic map, almost
anything can happen to the variety.  For example, any variety
can be desingularized via
a pullback by a proper map.  On the other hand, if we restrict ourselves to
equidimensional
germs of finite maps, it turns out that the geometry of the pullback can
``only get worse'' given the proper definition of ``worse.'' Specifically we
ask the question (note that every variety in this paper is complex analytic):

\begin{question} \label{theq}
Let $(W,0) \subset (\C^n,0)$ be a germ of a subvariety and $(F,0) \colon
(\C^n,0) \to (\C^n,0)$ be a germ of a finite holomorphic map.
We wish to find geometric
properties of $V := F^{-1}(W)$, such that if $(V,0)$ possesses property
$A$, then $(W,0)$ must possess property $A$.
\end{question}

The specific case with $A$ being the property of smoothness was asked
and partially solved by Ebenfelt and Rothschild~\cite{ER:finmap}.  Their
motivation was a question in CR geometry.  In particular they were interested
in the image $F(M)$ of a CR manifold $M$ for a finite map $F$.
It turns out a natural condition
to check is $\sM = \sF^{-1}(\sF(\sM))$, where $\sM$ is
the complexification of $M$, and $\sF$ is the
complexification of the finite map.  The natural question to
ask then is if $\sF(\sM)$ is a submanifold.

An obvious property answering the question is the property of irreducibility.
That is, if $V$ is irreducible as a germ at the origin, then $W$ must be
irreducible as a germ at the origin.

In this paper we will be looking at properties of the coordinate
ring $\sO_{V,0} = \sO_0 / I_0(V)$, where $\sO_0$ is the ring of germs of
holomorphic functions on $(\C^n,0)$ and $I_0(V) \subset \sO_0$ is the ideal
of germs of holomorphic functions vanishing on a subvariety $V$.
$\sO_{V,0}$ has no zero divisors
precisely when the germ $(V,0)$ is irreducible as a variety.
We say $V$ is \emph{normal} at $0$
if $\sO_{V,0}$ is integrally closed.  Equivalently, $V$ is \emph{normal} at $0$
if $\sO_{V,0} = \widetilde{\sO}_{V,0}$, where $\widetilde{\sO}_{V,0}$ is the
ring
of germs of weakly holomorphic functions on $V$ (essentially a function is
weakly
holomorphic if it is holomorphic on the smooth part of $V$ and locally
bounded).  Another interesting possible property of such rings is
factoriality (unique factorization property).  
We will say $V$ is \emph{factorial} (resp.\@ \emph{prefactorial})
at $0$ whenever $\sO_{V,0}$
is factorial (resp.\@ prefactorial).
A ring $R$ is said to be \emph{prefactorial} when every prime ideal of
height one is
a radical of a principal ideal.
In terms of varieties, this property
is equivalent to
every codimension one subvariety of $V$ being defined (as a set) by the
vanishing of a single function
in $\sO_{V,0}$.
In this paper, we prove that the two properties normal and prefactorial answer
Question~\ref{theq}.
We provide further conditions when a third property, being smooth, answers
the question.

\begin{thm} \label{mainthm}
Let $(W,0) \subset (\C^n,0)$ be a germ of a subvariety and 
$(F,0) \colon (\C^n,0) \to (\C^n,0)$ be a germ of a finite holomorphic map.
Let $V := F^{-1}(W)$.  If $(V,0)$ is normal (resp.\@ prefactorial)
then $(W,0)$ is normal (resp.\@ prefactorial).
\end{thm}

The ``normal'' part of the theorem was proved earlier by
Joseph Lipman (unpublished, see~\cite{ER:finmap}).
We present a different proof, using weakly
holomorphic functions, which is more digestible to
analysts.

Because of the simple structure of one-dimensional varieties the properties
of being factorial, normal, and smooth are equivalent.  On the other hand
every one-dimensional variety is prefactorial.
When the dimension is one we can therefore replace ``normal'' by ``smooth.''
There are also several other situations in which it can be proved that smooth
$V$ implies smooth $W$.  
We summarize results about the smooth case
proved by Ebenfelt and Rothschild~\cite{ER:finmap} 
together with results of the present paper
in the following theorem.
By $F$ being \emph{triangular} we mean that
the $k^{\text{th}}$ coordinate of $F$ depends only on the first $k$ variables.

\begin{thm} \label{smooththm}
Let $(W,0) \subset (\C^n,0)$ be a germ of a subvariety and 
$(F,0) \colon (\C^n,0) \to (\C^n,0)$ be a germ of a finite holomorphic map.
Let $V := F^{-1}(W)$ and suppose $(V,0)$ is smooth.  Suppose
one of the following conditions holds for some representatives $F$ and $V$:
\begin{enumerate}[(i)]
\item \label{st:onedim}
if $\dim(V) = 1$, or
\item \label{st:derv}
if $\det(DF)|_V \not\equiv 0$, or
\item \label{st:1to1}
if $F|_{V\setminus S}$ is one-to-one for some proper (possibly empty)
subvariety $S \subset V$, or
\item \label{st:prime}
$F$ is of multiplicity $p$ at $0$ ($F$ generically $p$ to one), where $p$ is
prime, or
\item \label{st:triang}
$F$ is triangular,
\end{enumerate}
then $(W,0)$ is smooth.
\end{thm}

Parts~\eqref{st:onedim} and~\eqref{st:derv} of the theorem
was proved in~\cite{ER:finmap}.  They 
proved~\eqref{st:onedim} directly using Puiseux parametrization.
There appears to be no consensus as to whether this theorem should hold
unconditionally.  Neither proof nor counterexample has been
found so far despite much effort.  Interestingly,
the theorem would be false if considered in
characteristic~$p$.  See \S~\ref{section:examples}.

Note that ``smooth $\Rightarrow$ normal and prefactorial,'' and
therefore we have some information in the case when $V$ is smooth.  In this
case we know that $W$ must be both normal and prefactorial.

The question can be stated purely algebraically.
Let $I := (g_1,\ldots,g_k)$ be an ideal in $\sO_0$.
Let $f = (f_1,\ldots,f_n)$ be a set of germs in $\sO_0$ such that the
corresponding ideal is of finite codimension ($f$ is a finite map).
Write
\begin{equation}
J := \sqrt{(g_1 \circ f, \ldots, g_k \circ f)} .
\end{equation}
Question~\ref{theq} can be reformulated as follows.
\emph{Find all properties such that if $\sO_0 / J$ has property $A$,
then $\sO_0 / I$ also has property $A$}.
The difficulty comes from taking the radical.  It is difficult to prove results
about the generators of the radical of an ideal even if we know much about
the ideal itself.

The author would like to acknowledge Peter Ebenfelt and
Linda Rothschild for suggesting the problem.  Also the author acknowledges
Jon Armel, John D'Angelo, Orest
Bucicovschi, J{\'a}ra Cimrman, and Xiaojun Huang, for helpful
discussions about this problem.
The author would also like
to acknowledge AIM where the special case of the smoothness was 
discussed during a workshop on complexity in CR geometry.


\section{Normal implies normal} \label{section:normal}

In this section we prove the normality part of Theorem~\ref{mainthm}.
We restate this claim for reader convenience.
 
\begin{thm} \label{normalitythm}
Let $(W,0) \subset (\C^n,0)$ be a germ of a subvariety and 
$(F,0) \colon (\C^n,0) \to (\C^n,0)$ be a germ of a finite holomorphic map.
Let $V := F^{-1}(W)$.  If $(V,0)$ is normal
then $(W,0)$ is normal.
\end{thm}

As mentioned before,
this theorem was suggested and proved by Joseph Lipman
using an algebraic method.
We present a different proof using analysis.
Following Whitney~\cite{Whitney:book}
we define $V$ to be normal at the origin if all weakly holomorphic functions
are in fact holomorphic.  This definition is equivalent to the more traditional
definition of normality saying that $\sO_{V,0}$ is integrally closed.
Let $V$ be of pure dimension $m$.  A
function $f$ is \emph{weakly holomorphic} if there exists
a proper subvariety $S \subset V$ of dimension strictly
less than $m$, such that
$f$ is a holomorphic function defined on $V \setminus S$
and such that $f$ is locally
bounded near every point of $V$.  We denote the
set of germs at $0$ of weakly holomorphic functions on $V$ by
$\widetilde{\sO}_{V,0}$.
Weakly holomorphic functions are in fact
restrictions to $V$ of meromorphic functions on the ambient space,
and 
$\widetilde{\sO}_{V,0}$ is the integral closure of 
$\sO_{V,0}$.

\begin{proof}[Proof of Theorem~\ref{normalitythm}]
Pick small neighbourhoods of the origin and representatives of $W$, $F$
and $V$ such that $V$ is normal at the origin and $V := F^{-1}(W)$.
Suppose $f \colon W \to \C$ is a
weakly holomorphic function defined on $W$.
Then $f \circ F$ is a weakly holomorphic function on $V$.  As $V$ is
normal, there exists a holomorphic function $g$ defined on a neighbourhood
of $V$ in $\C^n$ such that $g|_V = f \circ F$.  We may need to pick
a smaller representative $V$ to get $g$, and therefore we may have had to
pick a smaller representative $W$ and a representative $F$ defined for smaller
neighbourhoods.

Let $k$ be the multiplicity of $F$.  For an appropriate neighbourhood $U$ of
the origin define the function
$\varphi \colon U \subset (\C^n)^k \to \C$ by
\begin{equation}
\varphi (Z_1,\ldots,Z_k) = \frac{1}{k} (g(Z_1) + g(Z_2) + \cdots + g(Z_k)),
\end{equation}
where each $Z_j \in \C^n$.  We note that $\varphi$ is a symmetric
function and therefore, $\varphi \circ F^{-1}$ is a well defined
holomorphic function on a neighbourhood of the origin of $\C^n$.

When $z \in W$, then for all $p \in F^{-1}(z)$ we must have
$f(z) = g(p)$.  Therefore
$\left(\varphi \circ F^{-1}\right)|_W = f$.
Hence $f$ is a restriction of a holomorphic function defined on a
neighbourhood of $W$ and thus $W$ is normal at the origin.
\end{proof}


\section{Prefactorial implies prefactorial} \label{section:prefactorial}

We say that a variety $V$ is prefactorial at 0, if the ring
$\sO_{V,0}$ is prefactorial.  That is, every prime ideal of height one
is the radical of a principal ideal.  In other words, if every
codimension one subvariety $S \subset V$ is the intersection of $V$ and
a hypersurface.  If $V$ is irreducible,
a well known lemma (see~\cite{CRT}) says that if every
prime ideal of height one is principal (and not just a radical of a principal
ideal), then $\sO_{V,0}$ would in fact be factorial.

The most interesting case is when combining the property of being normal
and prefactorial.  An obvious question therefore is: are there normal but
nonprefactorial varieties?  To see an example take any codimension two
subvariety of $\C^n$ that is not a local set theoretic
complete intersection, but that is contained in a normal hypersurface.
The classical example~\cite{Whitney:book} suffices here.
Take the set $S$ defined
by
\begin{equation}
\rank
\begin{bmatrix}
u_1 & u_2 & u_3 \\
v_1 & v_2 & v_3
\end{bmatrix}
\leq 1 .
\end{equation}
$S$ is a codimension two subvariety of $\C^6$ and $S$ is not a set
theoretic complete intersection at the origin.  The hypersurface
$V$ defined by
$u_1 v_2 - u_2 v_1 = 0$ is normal and contains $S$.  Hence $V$ is not
prefactorial, but it is normal.

An example of a subvariety of $\C^3$
that is factorial at the origin (and hence normal and prefactorial)
but not smooth is defined by $x^2+y^3+z^5=0$.  See
for example Shafarevich~\cite{Shafarevich:bag1}.

%
%
%
%
%

\begin{thm}
Let $(W,0) \subset (\C^n,0)$ be a germ of a subvariety and 
$(F,0) \colon (\C^n,0) \to (\C^n,0)$ be a germ of a finite holomorphic map.
Let $V := F^{-1}(W)$.  If $(V,0)$ is prefactorial
then $(W,0)$ is prefactorial.
\end{thm}

\begin{proof}
Take small enough representatives of $W$, $V$, and $F$.
Let $S \subset W$ be a subvariety of codimension one in $W$ with $0 \in S$.
Then $F^{-1}(S)$ is a codimension one subvariety of $V$.  As $V$ is
prefactorial, there exists $f \in \sO_{V,0}$ such that
$F^{-1}(S)$ is precisely the zero set of $f$.  We extend
$f$ to be a function in $\sO_0$ and we call the extension again $f$.  Then
\begin{equation}
F^{-1}(S) = V \cap \{ f = 0 \} .
\end{equation}
$F$ is a finite map, and hence a proper map of a sufficiently small
neighbourhood of the origin onto another sufficiently small neighbourhood of
the origin.  By the proper mapping theorem, $F$ maps the
hypersurface $\{ f = 0 \}$ to a subvariety.  The image subvariety is
a hypersurface and thus there exists a function $g$
such that $F(\{f=0\}) = \{ g = 0 \}$.  Furthermore, we have
\begin{equation}
S = W \cap \{ g = 0 \} .
\end{equation}
We can consider $g$ to be a function in $\sO_{W,0}$.
\end{proof}


\section{Smoothness} \label{section:smoothness}

Let us restate the theorem dealing with the smoothness condition
(Theorem~\ref{smooththm}) for reader convenience.  We will prove the theorem
(except part \eqref{st:derv}) in this section.

\begin{thmnonum} 
Let $(W,0) \subset (\C^n,0)$ be a germ of a subvariety and 
$(F,0) \colon (\C^n,0) \to (\C^n,0)$ be a germ of a finite holomorphic map.
Let $V := F^{-1}(W)$ and suppose $(V,0)$ is smooth.  Suppose
one of the following conditions holds for some representatives $F$ and $V$:
\begin{enumerate}[(i)]
\item 
if $\dim(V) = 1$, or
\item 
if $\det(DF)|_V \not\equiv 0$, or
\item 
if $F|_{V\setminus S}$ is one-to-one for some proper (possibly empty)
subvariety $S \subset V$, or
\item 
$F$ is of multiplicity $p$ at $0$ ($F$ generically $p$ to one), where $p$ is
prime, or
\item 
$F$ is triangular,
\end{enumerate}
then $(W,0)$ is smooth.
\end{thmnonum}

As we said before, parts \eqref{st:onedim} and \eqref{st:derv} were
proved by Ebenfelt and Rothschild~\cite{ER:finmap}.  Our proof
of normality provides another proof of \eqref{st:onedim}, because for
one-dimensional varieties normality is equivalent to smoothness.
We start with a proof of \eqref{st:1to1}.

\begin{proof}[Proof of Theorem~\ref{smooththm} part \eqref{st:1to1}]
Suppose that $F|_{V\setminus S}$ is one-to-one.  As $V$ is a manifold it is
normal and hence we can apply Theorem~\ref{mainthm} to conclude that $W$
is normal.  The inverse $F^{-1}|_{W \setminus F(S)}$
is well defined and holomorphic outside
the singular set of $W$ (a one-to-one holomorphic function between complex
manifolds of same dimension is holomorphic).
Therefore $F^{-1}$ is 
is weakly holomorphic on $W$ and hence holomorphic by normality of $W$.
Thus $W$
and $V$ are biholomorphic and $W$ is a submanifold.
\end{proof}

We will need to define multiplicity and state a few standard results about it
before proving part \eqref{st:prime}.

Let $F \colon V \to W$ be a holomorphic map defined on an irreducible
$m$-dimensional
variety $V$ to a variety $W$ with $\dim W = m$.
Let $p \in V$ be a point.
Suppose there is a neighbourhood $U \subset V$ of $p$
such that $F^{-1}(F(p)) \cap \overline{U} = \{ p \}$.
We then define the multiplicity of $F$ at $p$ as
\begin{equation}
\mu_p(F) :=
\limsup_{q \to F(p)}
\# \left[ F^{-1}(q) \cap U \right] .
\end{equation}
We note that the multiplicity is independent of the choice of $U$ subject to
the condition above.  When $V$ is locally reducible, the multiplicity
is the sum of the multiplicities of $F$ restricted to the different
branches of $V$.
When the image is a complex manifold, then after perhaps taking a
smaller neighbourhood $U' \subset W$ of $F(p)$, $F$ is a $\mu_p(F)$-sheeted
analytic cover over $U'$.  In this case it can be proved that for each
integer $k$ the sets
\begin{equation} \label{levelsets}
\{ z \in U \mid \mu_z(F) \geq k \}
\end{equation}
are subvarieties of $U$.  Similarly if both the target and the image 
are normal, then the map is again an analytic cover and
the sets \eqref{levelsets} are again subvarieties.
See~\cite{Chirka:book} (pg.\@ 102) and~\cite{Stoll:mult}
for more information and
proofs of the above statements.

Usually, but not always,
when we restrict $F$ to a submanifold (or a subvariety),
the multiplicity must drop.  One
sufficient (but not necessary) condition is if the Jacobian is identically
zero on the submanifold.

\begin{lemma} \label{smlemma1}
Let $(F,0) \colon (\C^n,0) \to (\C^n,0)$ be a germ of a finite holomorphic map.
Let $(M,0) \subset (\C^n,0)$ be a germ of a complex submanifold.
Suppose that $\det (DF)|_M \equiv 0$, then $\mu_0(F|_M) < \mu_0(F)$.
\end{lemma}

\begin{proof}
Suppose that $k = \mu_0(F|_M) = \mu_0(F)$.  If we show that for a generic
point $p \in M$, $\mu_p(F) = 1$, then $F$ would have to be locally one-to-one
near generic points of $M$.
A one-to-one holomorphic map of complex submanifolds of same
dimension would be biholomorphic, which would contradict $\det(DF)|_M \equiv 0$.

We pick representatives of the germs such that
$F^{-1}(F(0)) = \{ 0 \}$, then for a generic point $p$ near 0,
$\# F^{-1}(F(p)) = k$.  Since
$\mu_0(F|_M)$ is also $k$, $\# F^{-1}(F(p)) = k$ for a generic point
$p$ on $M$
as well.  Take such a generic point $p \in M$,
and look at the sets
$F^{-1}(F(q))$ for points $q$ near $p$.
They all have cardinality $k$ and so
$F^{-1}$ is a continuous map to the symmetric space $(\C^n)^k_{sym}$
(a multifunction, see~\cite{Whitney:book}).
Then
there exists a neighbourhood $U \subset \C^n$
of $p$ such that $F^{-1}(F(q)) \cap U = \{ q \}$.  Hence $\mu_p(F) = 1$,
and the lemma is proved.
\end{proof}

\begin{lemma} \label{smlemma2}
Let $(F,0) \colon (\C^n,0) \to (\C^n,0)$ be a germ of a finite holomorphic map.
Let $(V,0)$ be an normal germ of a subvariety such that for some
representatives we have $V=F^{-1}(F(V))$.
Then $\mu_0(F|_V)$ divides $\mu_0(F)$.
\end{lemma}

The lemma does not hold in general without some extra assumptions on
$V$ such as normality and $V=F^{-1}(F(V))$.
Take the map $F$
defined by
$(z,w) \mapsto (z^aw^b,w^2-z^3)$.  The multiplicity of $F$ at the origin
is $2a+3b$ (see~\cite{Chirka:book} pg.\@ 109), but when restricted to the
manifold $\{ w=0 \}$, the map has multiplicity $3$.  When restricted to
$V:=F^{-1}(F(\{w=0\})) = \{ w = 0 \} \cup \{ z=0 \}$ (which is now reducible
and hence not normal),
the multiplicity is the sum of the multiplicities along both branches:
\begin{equation}
\mu_0(F|_V) =
\mu_0(F|_{\{z=0\}}) + 
\mu_0(F|_{\{w=0\}})
= 2 + 3 .
\end{equation}

\begin{proof}[Proof of Lemma~\ref{smlemma2}]
We take small representatives of $F$ and $V$ as before.
As we have noted, the sets $\{ z \mid \mu_z(F) \geq k \}$ are subvarieties
of some neighbourhood of the origin.  Hence, as $V$ is irreducible,
there exists a single integer $k$ and a proper subvariety $S \subset V$ (with
$\dim S < \dim V$) such
that $\mu_z(F) = k$ for all $z \in V \setminus S$.


Take a $p \in V \setminus S$ in this dense open subset such that
for $q \in F^{-1}(F(p))$ we also have $\mu_q(F) = k$ (which we can do since
$F^{-1}(F(p)) \subset V$ for all $p \in V$).  Let $z$ vary in a
small neighbourhood of $p$ in $\C^n$.  As $F^{-1}$ is a continuous map into
the symmetric space, each point $w \in F^{-1}(F(z))$ must split
$k$ times as $z$ (and hence each $w$ in $F^{-1}(F(z))$) leaves $V$,
because $\mu_w(F) = k$ when $w \in V$.
Thus
\begin{equation} \label{eq:split}
\Big( \# F^{-1}(F(p)) \Big) \times k = \mu_0(F) .
\end{equation}
Since $V$ is normal, hence $F(V)$ is normal by Theorem~\ref{mainthm}.
$F$ is locally an analytic cover.
Since \eqref{eq:split} holds for a generic $p \in V$, then as 
$F$ is an analytic cover we have
$\mu_0(F|_V) = \frac{\mu_0(F)}{k}.$  See Stoll~\cite{Stoll:mult}
Proposition 1.3.
%
%
\end{proof}

We can now prove part \eqref{st:prime} of the theorem.

\begin{proof}[Proof of Theorem~\ref{smooththm} part \eqref{st:prime}]
First note that if $\det (DF)|_V \not\equiv 0$, then the claim 
holds by part~\eqref{st:derv}.
In the opposite case, we apply Lemma~\ref{smlemma1} and Lemma~\ref{smlemma2}
to conclude that
$\mu_0(F|_M)$ divides $\mu_0(F)$ and
$\mu_0(F|_M) < \mu_0(F)$.  We assumed that $\mu_0(F)$ is prime and hence
$\mu_0(F|_M) = 1$.  Therefore $F$ is generically one-to-one on $M$
and we can apply part \eqref{st:1to1} of the theorem to finish the proof.
\end{proof}

Finally we prove part \eqref{st:triang}.  We note again what
we mean by $F$ being triangular.  We say $F = (F_1,F_2,\ldots,F_n)$
is triangular if $F_j$ depends only on $z_1,\ldots,z_j$.

\begin{proof}[Proof of Theorem~\ref{smooththm} part \eqref{st:triang}]
The proof follows by induction.  Suppose that \eqref{st:triang}
works for triangular maps in dimension $n-1$.  Let $F' =
(F_1,\ldots,F_{n-1})$ and $z' = (z_1,\ldots,z_{n-1})$.  We can factor
$F = G \circ H$, where
\begin{align}
H(z',z_n) &:= \big(z',F_n(z_n)\big) , \\
G(w',w_n) &:= \big(F'(w'),w_n\big) .
\end{align}
Since $F$ is finite, both $G$ and $H$ are finite.  Let $X = G^{-1}(W)$
and we note that $V = H^{-1}(X)$.

We claim that the conclusion of Theorem~\ref{smooththm} holds for
maps of the same form as $G$ and $H$.
That is, we first prove that $X$ must be smooth
and then that $W$ must be smooth.

We start with the map $H$.
If $V$
can be written as $z_n = f(z')$ then the image is $w_n = F_n(f(w'))$ and we are
done.  In all other cases, we consider smooth
one-dimensional curves $\Gamma$ through 0 in the
$z'$ space and restrict the map to $\Gamma \times \C$
to get a $\C^2$ to $\C^2$ map that is of the form
$(z, w) \mapsto (z, F_n(w))$.  Then $V \cap (\Gamma \times \C)$ is smooth and
we apply one-dimensional claim (part \eqref{st:onedim}) to note that
$W \cap (\Gamma \times \C)$ must be smooth.  As $\Gamma$ was arbitrary smooth
curve through 0 in the $z'$ space, 
$X$ is smooth at the origin.

To prove the claim for $G$, we apply the induction hypothesis by
first restricting to $\C^{n-1}$ by setting
$w_n = 0$.  The claim follows.
\end{proof}


\section{Examples} \label{section:examples}

If we take the philosophy that the geometry of the pullback $F^{-1}(W)$ must
be ``worse'' than that of $W$ itself, we should keep the following example in
mind.

\begin{example}
Take the map $F \colon \C^2 \to \C^2$ defined by
$(z,w) \mapsto (z^3+w,z^2)$.
Let $\zeta,\omega$ be the coordinates of the target.  Then the pullback of
the cusp $\omega^3 - \zeta^2 = 0$ is given by
\begin{equation}
0 = (z^2)^3 - (z^3+w)^2 = 2z^3w + w^2 = w(2z^3+w) .
\end{equation}
That is, the pullback is the intersection of two nonsingular lines.
\end{example}

As one would expect,
the line of reasoning outlined in this paper is not
easily applicable to real geometry.  

\begin{example}
Take the map $F \colon \R^2 \to \R^2$ defined by
$(x,y) \mapsto (x^3,x^2+y)$.
$F$ is one-to-one and the pullback of the cusp is the line $y=0$.  What goes
wrong here is, of course, the fact that for real $x$ the function $x^3$ is one
to one.
\end{example}

The following example by Melvin Hochster (see~\cite{ER:finmap})
means that the claim that smooth $V$ implies smooth $W$ will not be
shown by simple algebra only.

\begin{example}
Let us work in a field of characteristic $p$.  Let $x,y,z$ denote the 
coordinates on the source side and $u,v,w$ the coordinates on the
target.  Define $F$ by
\begin{equation}
u = x^p, \ \ v = xy+z, \ \ w = y^p .
\end{equation}
Then let $W = \{ v^p - uw = 0 \}$.
We compute the inverse image of $W$, using the fact that we are in
characteristic $p$,
\begin{equation}
0 = (xy+z)^p - x^py^p
= z^p .
\end{equation}
Therefore $V = F^{-1}(W) = \{ z = 0 \}$.
\end{example}


\section{The other direction} \label{section:otherdir}

We could also ask the ``inverse'' question.

\begin{question} \label{theqinv}
Let $(W,0) \subset (\C^n,0)$ be a germ of a subvariety and $(F,0) \colon
(\C^n,0) \to (\C^n,0)$ be a germ of a finite holomorphic map.
Let $V := F^{-1}(W)$.
We wish to find geometric properties of $(W,0)$,
such that if $(W,0)$ possesses property
$A$, then $(V,0)$ must possess property $A$.
\end{question}

It turns out very little geometric information of $W$ must be passed onto $V$
via the pull back.  The simplest property that is preserved is that of
reducibility.  The following example shows what can happen in general.

\begin{example}
Suppose that $V \subset \C^n$ is of codimension $k$ and is defined by $k$
holomorphic functions $f_1,\ldots,f_k$.  We can complete the defining
functions to a set $f_1,\ldots,f_k,f_{k+1},\ldots,f_n$ such that 0
is the only solution to $f_j(z) = 0$.  We let the mapping $F$
be defined by $F=(f_1,\ldots,f_n)$.  If we call $w$ the variables on the
target space then $F(V)$ is defined by $w_1 = \cdots = w_k = 0$.
Hence $W$ is a complex submanifold, and obviously $V = F^{-1}(W)$.
\end{example}

Any set theoretic complete intersection is the inverse image of a
submanifold.  Hence, essentially the only property which is preserved
under the inverse image by a finite holomorphic map is the property
of being a set theoretic complete intersection.


\def\MR#1{\relax\ifhmode\unskip\spacefactor3000 \space\fi%
  \href{http://www.ams.org/mathscinet-getitem?mr=#1}{MR#1}}

\end{document}